\documentclass[final]{siamonline1116}

\usepackage{lipsum}
\usepackage{amsfonts}
\usepackage{graphicx}
\usepackage{epstopdf}
\usepackage{algorithmic}
\ifpdf
  \DeclareGraphicsExtensions{.eps,.pdf,.png,.jpg}
\else
  \DeclareGraphicsExtensions{.eps}
\fi

\numberwithin{theorem}{section}

\newcommand{\TheTitle}{Function approximation using gradient information with application to parametric and stochastic differential equations}
\newcommand{\TheAuthors}{Gleb Ryzhakov and Ivan Oseledets}

\headers{\TheTitle}{\TheAuthors}

\title{{\TheTitle}%
}

\author{
  Gleb Ryzhakov\thanks{Skoltech, Moscow, Russia
    (\email{g.ryzhakov@skoltech.ru}).}
  \and
  Ivan Oseledets\thanks{Skoltech, Moscow, Russia
    (\email{i.oseledets@skoltech.ru}
    ).}
}

\usepackage{amsopn}
\usepackage{subcaption}
\usepackage[space,multidot]{grffile}
\usepackage{placeins}

\def\abs|#1|{\ensuremath{\left\lvert#1\right\rvert}}
\def\norm|#1|{\ensuremath{\left\lVert#1\right\rVert}}
\DeclareMathOperator{\grad}{grad}
\DeclareMathOperator{\rank}{rank}
\DeclareMathOperator{\mmod}{{ }mod}
\let\div=\relax
\DeclareMathOperator{\div}{div}
\def\bydef{\ensuremath{:=}}
\def\AG{\textsc{autograd}}

\graphicspath{{img/}}

\def\wpic{0.5}

\ifpdf
\hypersetup{
  pdftitle={\TheTitle},
  pdfauthor={\TheAuthors}
}
\fi

\begin{document}

\maketitle

\begin{abstract}
In the paper we consider the problem of multivariate function approximation in polynomial basis.
In order to solve this problem, we adjust the least squares method (LSM) by adding information about derivatives of the function.
This modification allows reducing the number of evaluations of approximating function while keeping the accuracy at the appropriate level.
We propose several techniques for time-efficient calculation of derivatives in various applications.
Numerical examples are given for comparison between the standard LSM and the proposed approach.
\end{abstract}

\begin{keywords}
  least squares, uncertainty quantification, function approximation, polynomial chaos expansion
\end{keywords}

\begin{AMS}
65D05,
65D15,
41A10,
35C11 
\end{AMS}
\section{Introduction}
 Multivariate
function approximation is a complex problem that arises in many practical applications.
We consider the
situation when the cost of a function calculation in a point is rather high (see~\cite{Z13,Z14} where polynomial approximation is applied to electric circuit).
One of the
standard
approaches to the approximation is
using
linear combination of elements from some polynomial basis as an approximant
and
the least squares method (LSM) to find the coefficients.
To reduce the complexity and
to improve the accuracy of the approximation it is needed to solve two problems,
how to decrease the number of sample points and how to appropriately choose them.

In our approach, the first problem is solved by calculation
of the function gradient.
The classical result by Baur and Strassen~\cite{BS83} on automatic differentiation states
that if the function is given as a composition of elementary operations (i.e.\@ it is a rational function),
then the gradient can be computed at most~$4N$ operations,
where~$N$ is the cost of a single evaluation of the function.
The method of computing derivatives is also known as ``backpropagation''.
This method allows fast calculating the gradient of the function that also consists of elementary functions and user-defined ones~\cite{G92}.
In our research we use Python package \AG~\cite{AG}
since it allows
calculating
derivatives in an efficient and convenient way.
However, in a case when the function is given as a black box, it is hard to find its derivatives with \AG.
In order to accelerate the calculation of derivatives in such cases, we use \emph{ad~hoc} methods.
These methods also enable calculating function in less number of points compared to the standard ones while keeping high accuracy.

The second problem of
sampling
can be solved in several ways.
We use the algorithm based on maximum volume concept~\cite{MO16,BM} in high-dimensional problems.
Additionally we use Latin Hypercube Sampling (LHS) (see~\cite{LHS}).
Monte-Carlo method was used as a reference.

Main contributions of the paper are
\begin{itemize}
\item we propose a new method of multivariate function approximation
with the use of
its derivatives;
\item
numerical experiments that demonstrate the excellence in the accuracy of the proposed approach over the standard one in a number of cases are conducted; 
\item we propose several techniques for
accelerated
calculation of derivatives
so
it can be reached
the
equivalent
accuracy in much less time
when the cost of function evaluation is high.
\end{itemize}
\section{Mathematical background}\label{sec:math}
\subsection{System matrix}
Consider the problem of a  multivariate polynomial approximation of a smooth function~$f(x)$,
$x=(x^{(1)},\, x^{(2)},\,\ldots,\, x^{(l)}) \in\mathbb R^l$.
One of the classical approaches is to use the linear combination of some linearly independent polynomials~$\{P_i\}$
\begin{equation}
    f(x)\approx\widehat f(x)=\sum_{i=1}^M\alpha_i P_i(x).
    \label{approx}
\end{equation}
Unknown coefficients~$\alpha_i$ in~\eqref{approx} can be found using the least squares method (LSM)
which minimize
the $l_2$ norm of the residual
in some selected set of points~$G=\{\xi_1,\,\xi_2,\,\ldots,\,\xi_m\}\subset \mathbb R^l$, $m\geq M$
\begin{equation*}
    \sum_{i=1}^m\left( f(\xi_i) - \widehat f(\xi_i) \right)^2 \to\min_{\{\alpha_i\}}.
\end{equation*}

Let the values of the function~$f$ and its gradient be known
at the points of the set~$G$.
We use these values and the LSM approach to solve the following problem
\begin{equation}
    A \alpha = F,
    \label{main_sys}
\end{equation}
where the vector~$\alpha=\{\alpha_1,\, \alpha_2,\, \ldots,\, \alpha_m\}$ is the vector of unknown coefficients of decomposition~\eqref{approx};
$F\in\mathbb R^{(l+1)m}$ is the vector of values of the function~$f$ and its gradient at the points of the set~$G$
\begin{equation}
    F_i=
    \left\{
        \begin{aligned}
            f(\xi_i)           &, & \;\; 0 < i&\leq m,\\
            \partial_1 f(\xi_{i \mmod m})&, & m < i&\leq 2m\\
            \partial_2 f(\xi_{i \mmod m})&, & 2m < i&\leq 3m\\
            \cdots\\
            \partial_l f(\xi_{i \mmod m})&, & l m < i&\leq (l+1)m,
        \end{aligned}
    \right.
    \label{func_and_derive}
\end{equation}
here~$\partial_i$ is~$i^{\text{th}}$ gradient's component
\begin{equation*}
    \partial_i\bydef\frac{\partial}{\partial x^{(i)}}.
\end{equation*}
The matrix~$A=\{a_{ij}\}$ has the following elements
\begin{equation}
    a_{ij}=
    \left\{
        \begin{aligned}
            P_j(\xi_i)           &, & \;\;0 < i&\leq m,\\
            \partial_1 P_j(\xi_{i \mmod m})&, & m < i&\leq 2m\\
            \partial_2 P_j(\xi_{i \mmod m})&, & 2m < i&\leq 3m\\
            \cdots\\
            \partial_l P_j(\xi_{i \mmod m})&, & l m < i&\leq (l+1)m
        \end{aligned}
    \right..
    \label{system:mat}
\end{equation}
In this approach we use additional information about the function (derivatives).
Thus the resulting matrix has more rows than without using derivatives.
So we expect that
adding a new column to the system matrix does not change its full-column rank property.
That means
that we can use polynomials of a higher degree to increase the accuracy of the approximation.

As a set of polynomials~$\{P_i\}$ we take a subset of the tensor product of
unidimensional
polynomials~$f_i$ of the form
\begin{equation*}
    f_{\beta_1}(x^{(1)})
    f_{\beta_2}(x^{(2)})
    \cdots
    f_{\beta_l}(x^{(l)})
\end{equation*}
with hyperbolic transaction scheme (see review~\cite{DTU16}) where polynomials degree~$\{\beta_i\}$ satisfy inequality
\begin{equation*}
    \sum_{i=1}^l\beta_i^p\leq q^p, \qquad \beta_i\in\mathbb N\cup\{0\},
\end{equation*}
for some~$p\in(0,1]$.
In the case~$p=1$, total number of elements of
the
polynomial set is equal to~$(l+p)!/(l!p!)$.

\section{Algorithm of multivariate function  approximation using derivatives}
We propose the following algorithm for the method
of function approximation
with the use of function derivatives based on the  ideas described in Sec.~\ref{sec:math}.

\begin{algorithm}[h!]
    \caption{Algorithm for function approximation using the information of its derivatives}
\label{alg:main}
\renewcommand{\algorithmicrequire}{\textbf{Input:}}
\renewcommand{\algorithmicensure}{\textbf{Output:}}
\begin{algorithmic}[1]
    \REQUIRE{Number of points~$m$, function evaluator, function derivatives evaluator, function domain}
    \ENSURE{Coefficients of the expansion~$\alpha$}
\STATE Select $N\geq m$ nodes $\{\xi_1,\,\xi_2,\,\ldots,\,\xi_N\}$ from the function domain using one of the classical algorithms
\STATE If $N>m$, use \textsc{maxvol} algorithm to reduce the number of the selected points to~$m$ and select some subset of the set from the first step
\STATE Compute right-hand side~$F$ using~\eqref{func_and_derive}
\STATE Form the system matrix~$A$ using~\eqref{system:mat}
\STATE Solve the least squares problem $\alpha=A^+F$
\end{algorithmic}
\end{algorithm}

At the first step of the approximation of the given function~$f$
we choose the points where the function is evaluated.
We use random sampling, uniform distribution, Latin Hypercube Sampling (LHS) in different experiments.
In some cases we reduce the number of selected points using \textsc{maxvol}~\cite{BM} procedure.

Then we form a vector of the right-hand side of the overdetermined system~\eqref{main_sys} by calculating function values and its derivatives at the selected points.
The choice of a method that we used for derivative calculation depends on the task.
Different approaches are described below in Sec.~\ref{sec:ODE}--\ref{sec:PDE}.

The resulting coefficients~$\alpha$ can be used to evaluate the values of the function approximant directly.
They can also be used
to find some properties of the unknown function such as the mean value and standard deviation if the function is random one as in the example in Sec.~\ref{sec:noise}.

The complexity of the algorithm depends on the complexity of function and
calculation of its derivatives.
We assume that an asymptotic complexity of derivatives calculation is the same as the function calculation one,
i.e.\@ the function satisfies the Baur and Strassen theorem~\cite{BS83}.
Theoretically we need no more then~$4N$ operations to find all derivatives if the function is rational one, where~$N$ is the number of operations necessary to find  the function.
Several effective ways for fast computation of the derivatives are given below in the following sections.
If the function evaluation at one point takes a long time, our approach allows to dramatically
decrease
the total calculation time in the case of high-dimensional tasks ($d\gg4$) to get the same size of the system matrix~$A$.
We do not know the exact estimates for the Lebesgue constant for our approach,
so we can not theoretically estimate the accuracy of the approximation.

Now we demonstrate the quality of the proposed approach on model examples.

\section{Model example of approximation}
We consider the following function of two variables
\begin{equation*}
    f(x_1,\,x_2)=\exp\left(-x_1^2-0.5 (x_2-1) x_2\right)
\end{equation*}
to check interpolation error with and without using the information of its derivatives.

The interpolation is performed on the rectangle~$[-2,\,2]^2$.
We use different number of points $N=50$, $200$ where the function values  and its gradient are calculated.
We use~\textsc{maxvol} algorithm to select the desired amount of points from the initial set of~$10^4$ points uniformly distributed on the rectangle.
Chebyshev polynomials are used as the basis functions.

Results are shown in Fig.~\ref{Fig:interp1}--\ref{Fig:interp2},
where $L_2$, relative error of the interpolation, is given as a function of the number of the monomials.
\def\cftw{0.45}%
\newcounter{figcount}%
\setcounter{figcount}{0}%
\def\figsimpl#1#2#3{%
\refstepcounter{figcount}%
\begin{figure}[tb]
\centering
\begin{subfigure}[h]{\cftw\textwidth}
    \includegraphics[width=\textwidth]{Num_end._pnts_#1,_Monom_type_cheb,_Range_0.5#2.pdf}
\caption{Range $[-0.5, \,0.5]$}
\end{subfigure}
\begin{subfigure}[h]{\cftw\textwidth}
    \includegraphics[width=\textwidth]{Num end. pnts=#1, Monom type=cheb, Range=1.0#2.pdf}
    \caption{Range $[-1,\,1]$}
\end{subfigure}
\begin{subfigure}[h]{\cftw\textwidth}
    \includegraphics[width=\textwidth]{Num end. pnts=#1, Monom type=cheb, Range=2.0#2.pdf}
    \caption{Range $[-2,\,2]$}
\end{subfigure}
\caption{Number of points $N=50$, #3 noise}
\label{Fig:interp\thefigcount}
\end{figure}%
}%
\def\figsimpl#1#2{%
\refstepcounter{figcount}%
\begin{figure}[tb]
\centering
\begin{subfigure}[h]{\cftw\textwidth}
    \includegraphics[width=\textwidth]{Num_end._pnts_50,_Monom_type_cheb,_Range_2.0#1.pdf}
\caption{Number of points is $N=50$}
\end{subfigure}
\begin{subfigure}[h]{\cftw\textwidth}
    \includegraphics[width=\textwidth]{Num_end._pnts_200,_Monom_type_cheb,_Range_2.0#1.pdf}
\caption{Number of points is $N=200$}
\end{subfigure}
\caption{Error of interpolation #2 noise}
\label{Fig:interp\thefigcount}
\end{figure}%
}%
\figsimpl{}{without}%
\figsimpl{_n}{with}%
One can see that the derivative approach gives more accurate results.

Now we will show how this approach can be applied to parametric ODEs.
\FloatBarrier
\section{Parametric Ordinary Differential Equations}\label{sec:ODE}
\subsection{Mathematical statement}
Consider a linear differential algebraic equation (DAE)
with coefficients depending on~$l$ parameters $\vec \xi=[\xi_1,\,\xi_2,\,\ldots,\,\xi_l]$
\begin{equation}
    E(\vec\xi)\frac{d\vec x(t,\vec\xi)}{dt}+A(\vec\xi)\vec x(t,\vec\xi)=B \vec u(t).
    \label{main_eq}
\end{equation}
Here~$E$ and~$A$ are $n\times n$ matrices,
$\vec x(t,\vec\xi)\in\mathbb R^n$ is the vector of unknowns,
$\rank E\leq n$,
$\rank A=n$.
The matrix~$B$ is $n\times d$,
$\vec u(t)\in\mathbb R^d$ is the vector of inputs.
This kind of equations describes, for example, electric circuits which contain only linear elements.
In the case of circuits with non-linear elements, we first linearise it on each time step and still get~\eqref{main_eq}
where matrices will be time-dependent.

As an example
some input parameters of the equation~\eqref{main_eq} can be random values.
The problem is to find the distribution of the solution.
Physically it means that some circuit parameters are not known exactly because of the nature of the environment,
and we want to know some properties of the output signal.
See~\cite{Z13} for examples of the polynomial chaos approach to several electric schemes.
Now we will show how to obtain the derivatives of the solution $\vec x$ efficiently.
\subsection{Obtaining derivatives}
\label{SS:derives}
\subsubsection{Derivatives via new right-hand side of the same equation}
Differentiating the equation~\eqref{main_eq} with respect to the $i^\text{th}$ parameter yields
\begin{equation}
    E(\vec\xi)\frac{d\vec x_i(t,\vec\xi)}{dt}+A(\vec\xi)\vec x_i(t,\vec\xi)= 
    B_i \vec u(t)
    -E_i(\vec\xi)\frac{d\vec x(t,\vec\xi)}{dt}
    -A_i(\vec\xi)\vec x(t,\vec\xi),
    \label{diff_eq}
\end{equation}
where
\begin{equation*}
    \mathcal X_i := \frac{\partial }{\partial\xi_i}\mathcal X \;\text{ for }\;
    \mathcal X = x(t, \vec\xi), E(\vec\xi), A(\vec\xi), B(\vec\xi).
\end{equation*}
So, if we have solved the equation~\eqref{main_eq}, then we can solve~\eqref{diff_eq} with a known right-hand side
to obtain the derivative with respect to the parameters.

\subsubsection{Derivatives via solving the extended system}
Consider the following equation
\begin{equation}
    \widetilde E(\vec\xi)\frac{d\vec x(t,\vec\xi)}{dt}+\widetilde A(\vec\xi)\widetilde{\vec x}(t,\vec\xi)=\widetilde B u(t),
    \label{eq-diff}
\end{equation}
with matrices
\begin{equation}
    \widetilde E =
    \begin{pmatrix}
        E   & 0 & 0 & \dots & 0\\
        E_1 & E & 0 & \dots & 0\\
        E_2 & 0 & E & \dots & 0\\
        \hdotsfor{5} \\ 
        E_l & 0 & 0 & \dots & E
    \end{pmatrix}
    ,\quad
    \widetilde A =
    \begin{pmatrix}
        A   & 0 & 0 & \dots & 0\\
        A_1 & A & 0 & \dots & 0\\
        A_2 & 0 & A & \dots & 0\\
        \hdotsfor{5} \\ 
        A_l & 0 & 0 & \dots & A
    \end{pmatrix}
    ,\quad
    \widetilde B =
    \begin{pmatrix}
        B  \\
        0 \\
        0 \\
        \cdots \\
        0\\
    \end{pmatrix}.
    \label{big_mats}
\end{equation}
Then the solution of this equation will be the vector
\begin{equation*}
    \widetilde{\vec x} = \{\vec x;\,\vec x_1;\, \vec x_2;\, \ldots;\, \vec x_l \}.
\end{equation*}

We can easily obtain the LU~factorization of the matrix with such a block structure as in~\eqref{big_mats}.
Indeed, denote the LU~factorization of the block~$A$ as $A=LU$. Then $\widetilde A=\widetilde L \widetilde U$ and
\begin{equation}
    \widetilde U =
    \begin{pmatrix}
        U   & 0 & 0 & \dots & 0\\
        0 & U & 0 & \dots & 0\\
        0 & 0 & U & \dots & 0\\
        \hdotsfor{5} \\ 
        0 & 0 & 0 & \dots & U
    \end{pmatrix}
    ,\quad 
    \widetilde L =
    \begin{pmatrix}
        L   & 0 & 0 & \dots & 0\\
        L_1 & 0 & 0 & \dots & 0\\
        L_2 & 0 & 0 & \dots & 0\\
        \hdotsfor{5} \\ 
        L_l & 0 & 0 & \dots & 0
    \end{pmatrix},
    \label{LU-sparse}
\end{equation}
where matrices~$L_i$ for $i=1,\dots, l$ are such that
\begin{equation*}
    A_i = L_i U, \quad i=1, \dots, l.
\end{equation*}
and they can be calculated fast since~$U$ is an upper triangular matrix and each~$A_i$ is sparse.

Let the matrix~$A$ be of size~$n\times n$,
number of operations required for solving system with the matrix~$\widetilde A$  is
\begin{equation}
    N = O(n^3) + O(c_1 n^2) + O(l n^2) + O(c_2).
    \label{num_of_oper}
\end{equation}
Here~$c_1$ is the total number of non-zero rows in the matrices~$A_i$,
$c_2$ is the total number of non-zeros elements in matrices~$A_i$.
The first two terms in~\eqref{num_of_oper} correspond to the LU-decomposition of the matrix~$\widetilde A$,
the last two terms are the number of operations required for the solution step.

\subsection{Some numerical examples}
In this section we present several examples of fast calculation of the function derivatives using the approaches described above.
Their timings are given as well as the asymptotic complexity.
\subsubsection{Linear electrical scheme}
\begin{figure}[h]
    \centering
    \includegraphics[width=\linewidth]{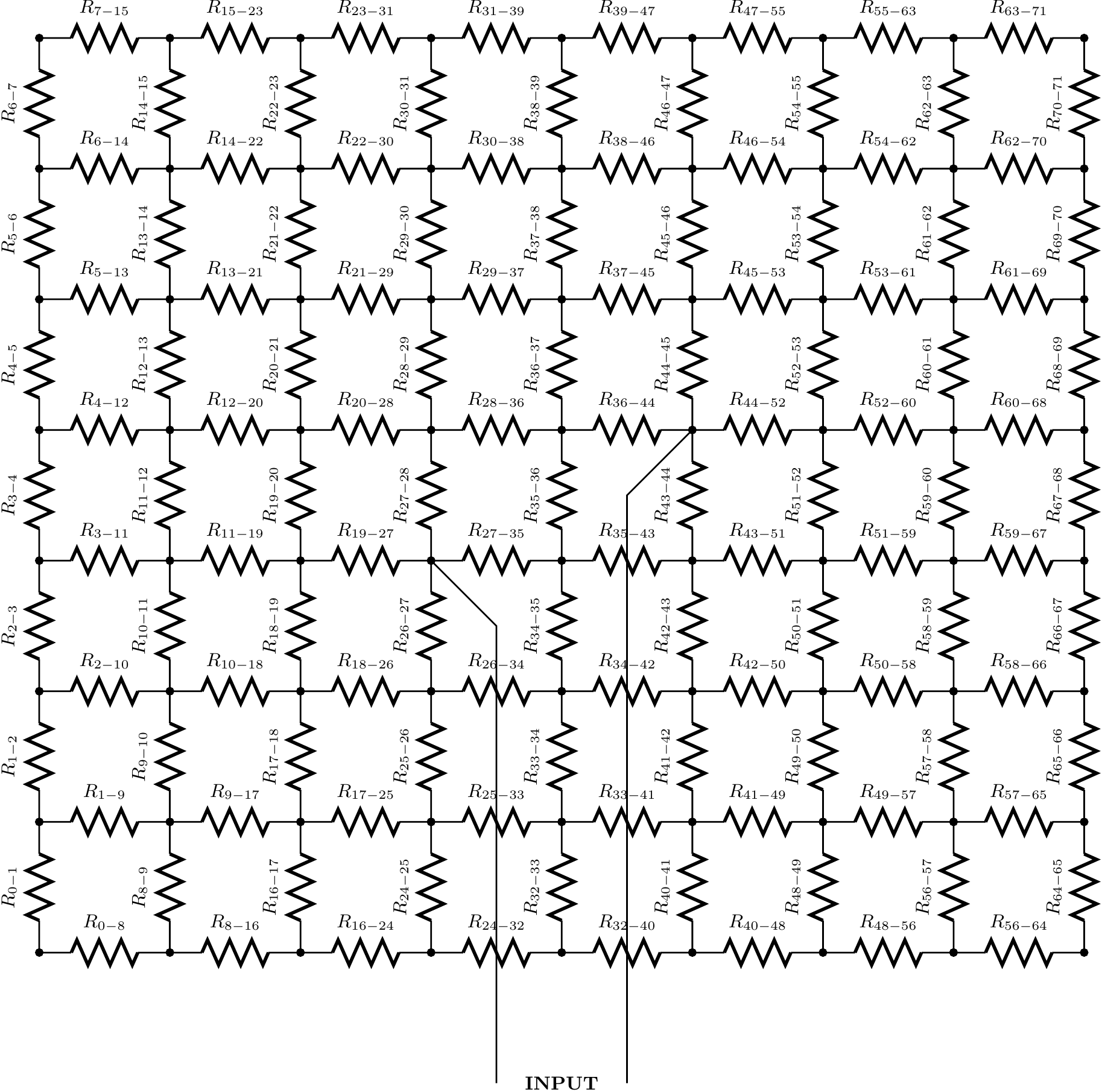}
    \caption{Testing scheme.}
    \label{fig:scheme1}
\end{figure}

Consider the linear scheme in Fig.~\ref{fig:scheme1} as an example.
The electrical schemes that we use in the experiments consists of
$967$ and $7687$
resistors,
only central part of the scheme is shown in Fig.~\ref{fig:scheme1}.
The closest elements to the inputs points are treated as unknown parameters in the sense that we take derivatives with respect to them.
For example, the first ten elements are listed below
\begin{gather*}
R_{35-36},\, R_{27-35},\, R_{28-36},\, R_{35-43},\, R_{36-44}, \\
R_{27-28},\,
R_{34-35},\,
R_{36-37},\,
R_{43-44},\,
R_{26-27}.
\end{gather*}
\begin{table}%
\centering%
\begin{tabular}[ht]{|p{4.5em}|c|c|c|c|c|}\hline
    Number of parameters~$l$ & 5    & 8    & 12   &  300  & 600  \\\hline
    Time (sec.)              & 0.07 & 0.11 & 0.16 &  3.96 & 8.18 \\\hline
\end{tabular}%
\caption{Built-in sparse LU-solver timings}%
\label{tbl:eq-diff}%
\end{table}%
\begin{table}%
\centering%
\begin{tabular}[ht]{|p{7em}|c|c|c|c|}\hline
    Number of parameters~$l$ & 2    & 10    & 300 &  600  \\\hline
    Time (sec. $\times 10^8$) for number of elements $k=967$         & 0.40 & 0.47 &  1.41 & 2.07  \\\hline
    Time (sec. $\times 10^{10}$) for number of elements $k=7687$     & 1.99 & 2.01 &  2.66 & 3.30  \\\hline
\end{tabular}%
\caption{Autograd timings}%
\label{tbl:autograd}%
\end{table}%
The timings of solving the equation~\eqref{eq-diff} with a different number of unknown parameters
are shown in Table~\ref{tbl:eq-diff}.
Matrices were stored in
Compressed Sparse Row (CSR)
format,
and the solution was obtained with Python build-in LU-factorization.
We can see that the
times grows almost linear  in the number
of the parameters when the block structure of the matrices is not accounted.
The number of operations for the solution of the system with the LU-decomposition~\eqref{LU-sparse}
of the matrix~$\widehat A$
for the electric circuit in Fig.~\ref{fig:scheme1} is shown in Table~\ref{tbl:oper}.%
\begin{table}%
\centering%
\begin{tabular}[ht]{|p{7em}|c|c|c|c|}\hline
    Number of parameters~$l$ & 2    & 10   & 300  & 600  \\\hline
    Time (sec.)          & 1415 & 1437 & 1366 & 1332 \\\hline
\end{tabular}
\caption{Timings for sparse system matrix solving}%
\label{tbl:oper}%
\end{table}%

We used the scheme with~$647$ elements to measure the time in the case of calculating derivatives using \AG.
The timings are shown in Table~\ref{tbl:autograd}.
We can see that the timings are practically independent of the number of arguments.

\subsubsection{Non-liner electrical scheme with noise}\label{sec:noise}
\begin{figure}[!th]
\centering
\includegraphics[width=\wpic\linewidth]{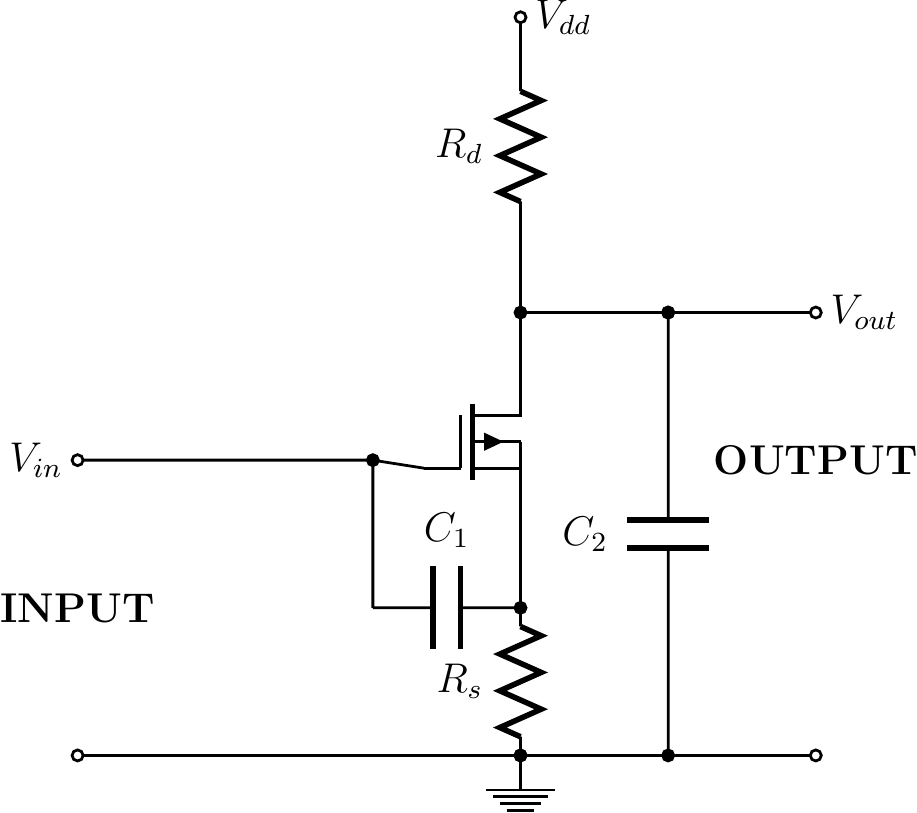}
\caption{Common-source amplifier}
\label{fig:scheme}
\end{figure}

We can apply our approach to model schemes that also have non-linear elements.
Consider a common-source amplifier shown in Fig.~\ref{fig:scheme}.
We use the MOSFET model of transistor in our numerical simulation.
Python package \textsc{ahkab}~\cite{AHKAB} was used for this purpose.
In the scheme~\ref{fig:scheme} the input voltage~$V_{\text{in}}$ is constant,
and the voltage~$V_{\text{dd}}$ is equal to~$3.3$ volts.
The values of parameters~$R_{\text{d}}$ and~$R_{\text{s}}$  have normal distribution
\begin{equation*}
    R_{\text{d}}\sim\mathcal N(10000,\,1000),\quad
    R_{\text{s}}\sim\mathcal N(300,\,150).
\end{equation*}
The temperature of the transistor is also random, $T\sim\mathcal N(300,\,70)$.
\def\cftw{0.55}%
\begin{figure}[!ht]
\centering
\begin{subfigure}[h]{\cftw\textwidth}
    \includegraphics[width=\textwidth]{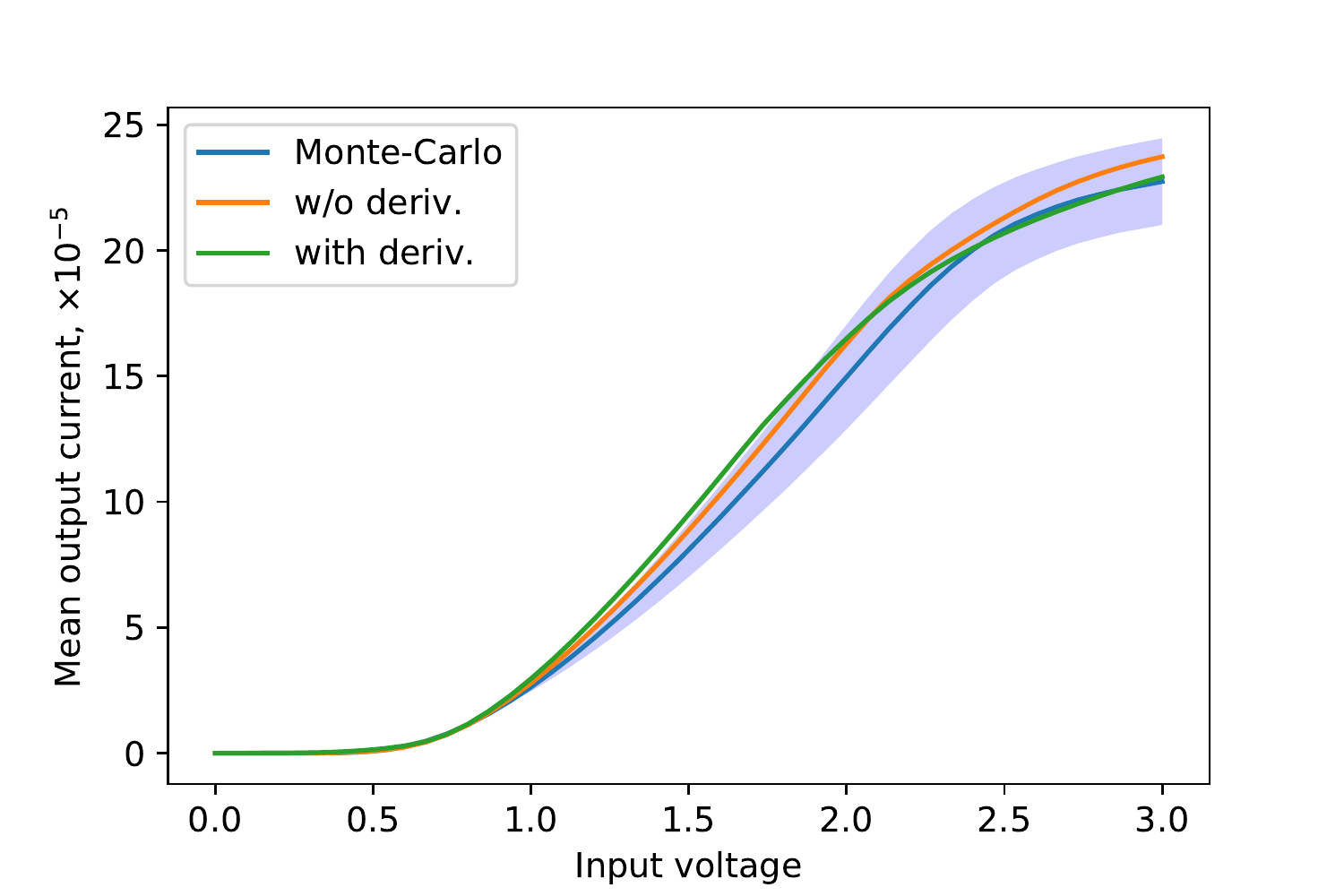}
    \caption{Mean value of output voltage as a function of input voltage.}
\end{subfigure}
\begin{subfigure}[h]{\cftw\textwidth}
    \includegraphics[width=\textwidth]{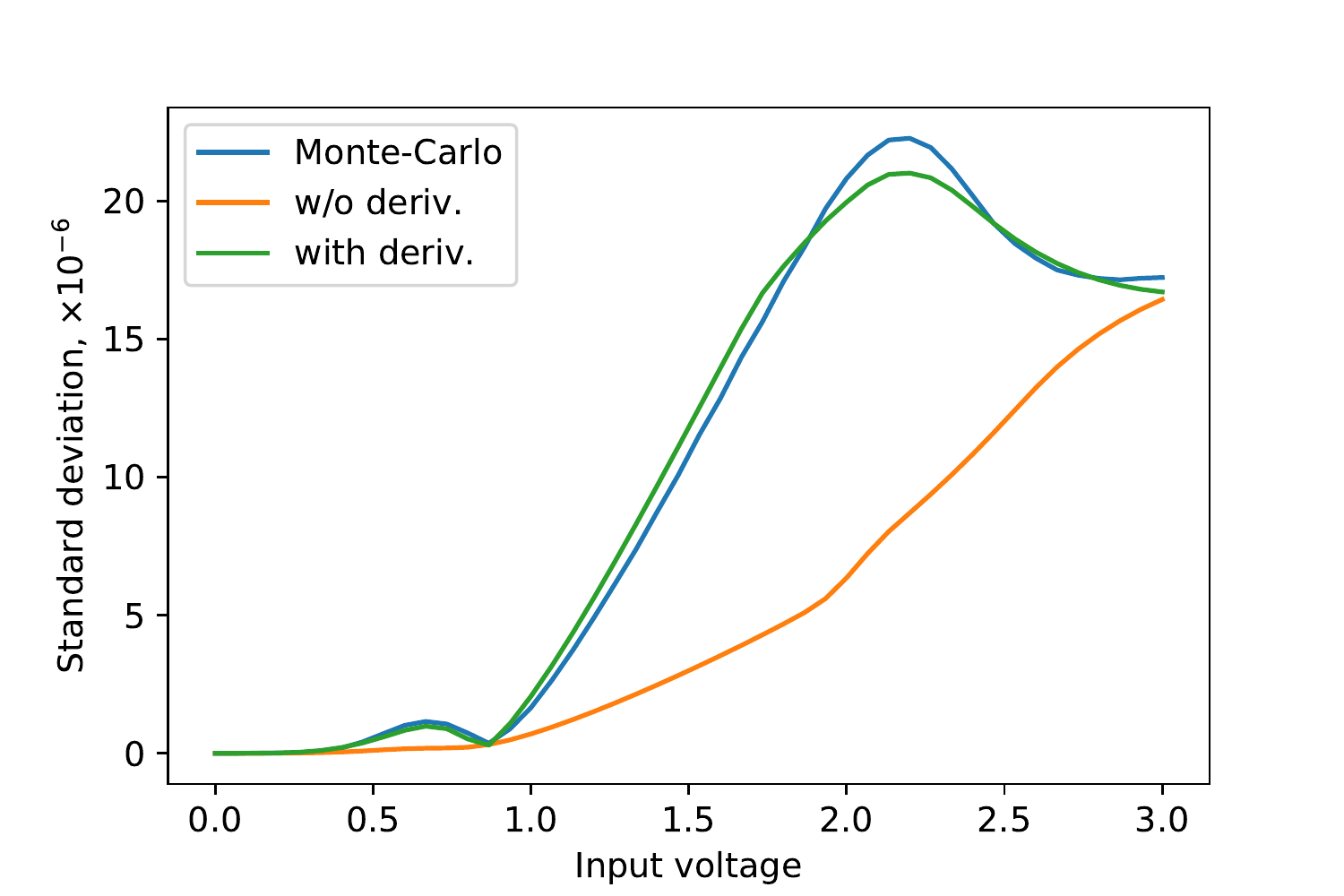}
    \caption{Standard deviation of output voltage as a function of input voltage.}
\end{subfigure}
\caption{Results of CS simulation.}
\label{fig:CS}
\end{figure}
The value of output voltage is the function of interest.
We approximate it by polynomials using Algorithm~\ref{alg:main}.
The technique described is Sec.~\ref{SS:derives} is used to obtain derivatives.
We use tensor product of univariate Hermite polynomials with weight function $w(x)=\exp\left(-x^2/2\right)$
as a set of basis functions
because input random parameters have a normal distribution and their
probability density function has the same form as $w(x)$.
The mean value of the output signal coincides with the first component of~$\alpha$,
i.e.\@ with the coefficient of the first polynomial~$P_0(x)=0$.
The standard deviation is equal to the square root of the sum of squares of
the remaining components
of the vector~$\alpha$
\begin{equation*}
\text{s.d.} =
\left(
\sum_{i=2}^M\alpha_i^2
\right)^{1/2}
\end{equation*}

We compare the results of simulation in  
in two cases: with and without using derivatives.
The number of points is equal to~$9$ in both cases.
In the case without using derivatives, we take~$3$ basis polynomials, so the system matrix has size~$9\times3$.
In the case of using derivatives, the number of basis polynomials is equal to~$5$, the size of the matrix is~$32\times5$.
Monte-Carlo simulation with $10^5$ samples is used as a reference.
The results of the model simulation are shown in Fig~\ref{fig:CS}.
The mean values of the output voltage are close to each other,
but one can see the advantages of the derivative approach
when the values of derivatives are used along with the values of the function.
Such small values of number of points and number of basis polynomials are used
to show the difference between methods.

\section{Partial differential equation}\label{sec:PDE}
Our method can also be applied to the classical uncertainty quantification UQ problem for the partial differential equation.
Let us consider the following differential equation, describing the stationary process of heat transfer
\begin{equation*}
    \div \left(k\,\grad u\right) = -F,
\end{equation*}
where the function~$k$ is a given thermal conductivity,
the function~$F$ is a given heat source
and~$u$ is an unknown temperature distribution.

We consider the case when~$k$ represents a lognormal random field
\begin{equation}
    k(x)=\exp(a+bg(x))
    \label{field_k}
\end{equation}
where~$g$ is a standard Gaussian random field with
the following correlation function
\begin{equation}
    \rho(x,\,y)=\exp\left(-\norm|x-y|^2/\sigma^2\right).
    \label{field_corr}
\end{equation}

For numerical experiments,
the values of parameters~$a$ and~$b$ in~\eqref{field_k} are $a=-0.0430888$, $b=0.29356$,
so the mean value and the standard deviation of~$k$ are~$1$ and~$0.3$, respectively.
In~\eqref{field_corr}, $\sigma=0.2$.

To model a random distribution of the thermal conductivity~$k$, we use
the Expansion Optimal Linear Estimation method~\cite{EOLE}.
According to this method, a random field is approximated by a truncated sum of the following form
\begin{equation}
    g(x) = \sum_{i=1}^N\frac{\xi_i}{\sqrt{l_i}}\left(\phi_i,\,C_{x\eta}\right)
    \label{field_expansion}
\end{equation}
where i.i.d.\@ $\xi_i\sim\mathcal N(0,1)$,
$C_{x\eta}$ are vectors with elements~$C^i_{x\eta}=\rho(x,\,\eta_i)$,
$l_i$ and $\phi_i$ are the eigenvalues and eigenvectors of the matrix~$C_{\eta\eta}$ accordingly.

The number of terms in~\eqref{field_expansion} is~$N=40$ for the numerical experiments.
The objective function~$T$ is the average temperature in the circle~$C$
with the center~$(-0.5,\,-0.5)$ and the radius~$0.2$
\begin{equation*}
    T=\int_C u\,ds,\quad
    C=\{x\in\mathbb R^2 \colon \left(x+0.5\right)^2 + \left(y+0.5\right)^2 < 0.2^2\}.
\end{equation*}
The heat source~$F$ is equal~$500$ in the circle with the center~$(0.5,\,0.5)$ and the radius~$0.2$,
and~$F=0$ elsewhere.
The Python package \textsc{dolfin-adjoint}~\cite{DAdj} is used for solving the partial differential equation and to obtain derivatives.

\begin{figure}[htb]%
\centering%
\includegraphics[width=\linewidth]{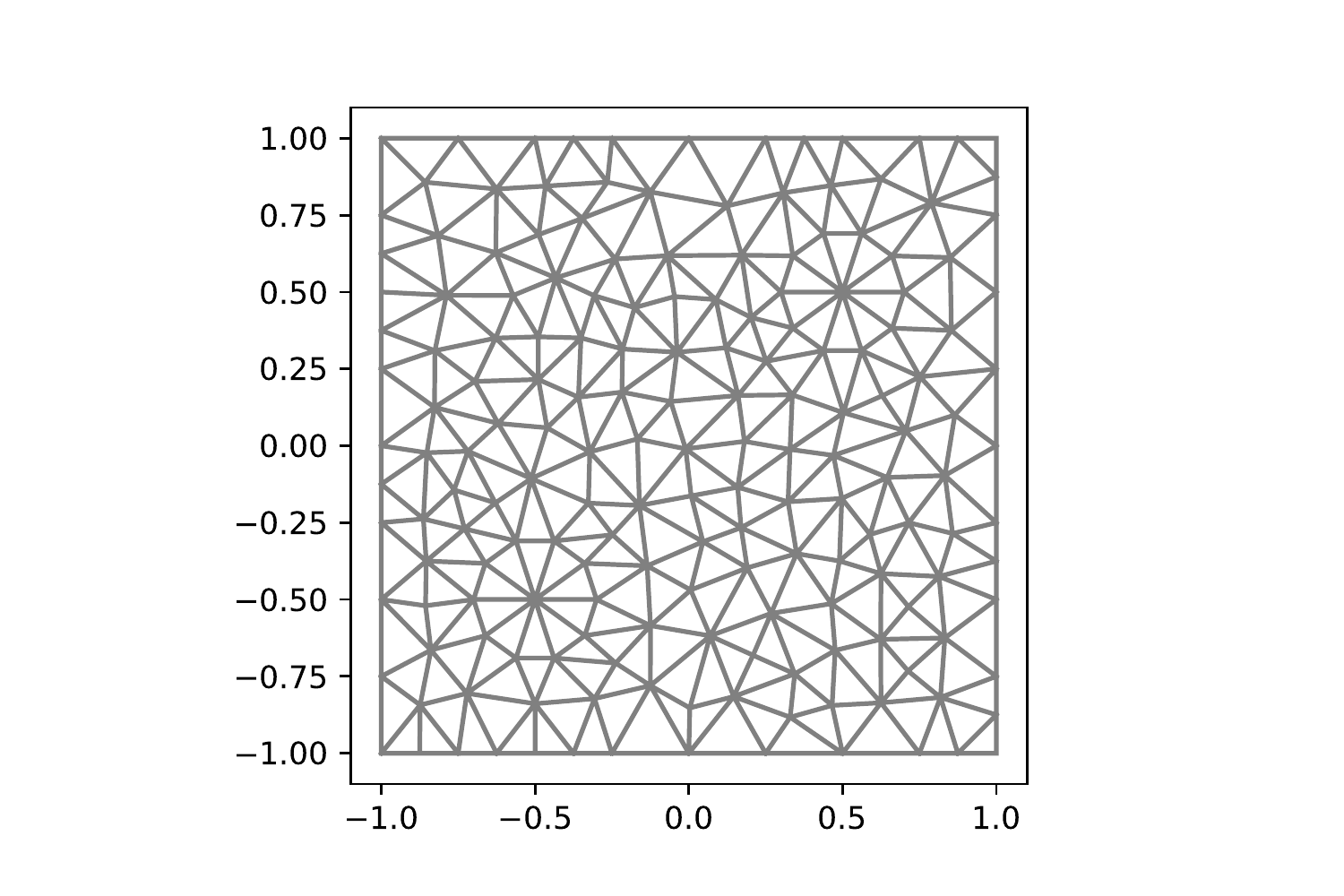}%
\caption{Mesh for finite-element method.}%
\label{fig:RF_mesh}%
\end{figure}

\begin{figure}[htb]%
\centering%
\includegraphics[width=\linewidth]{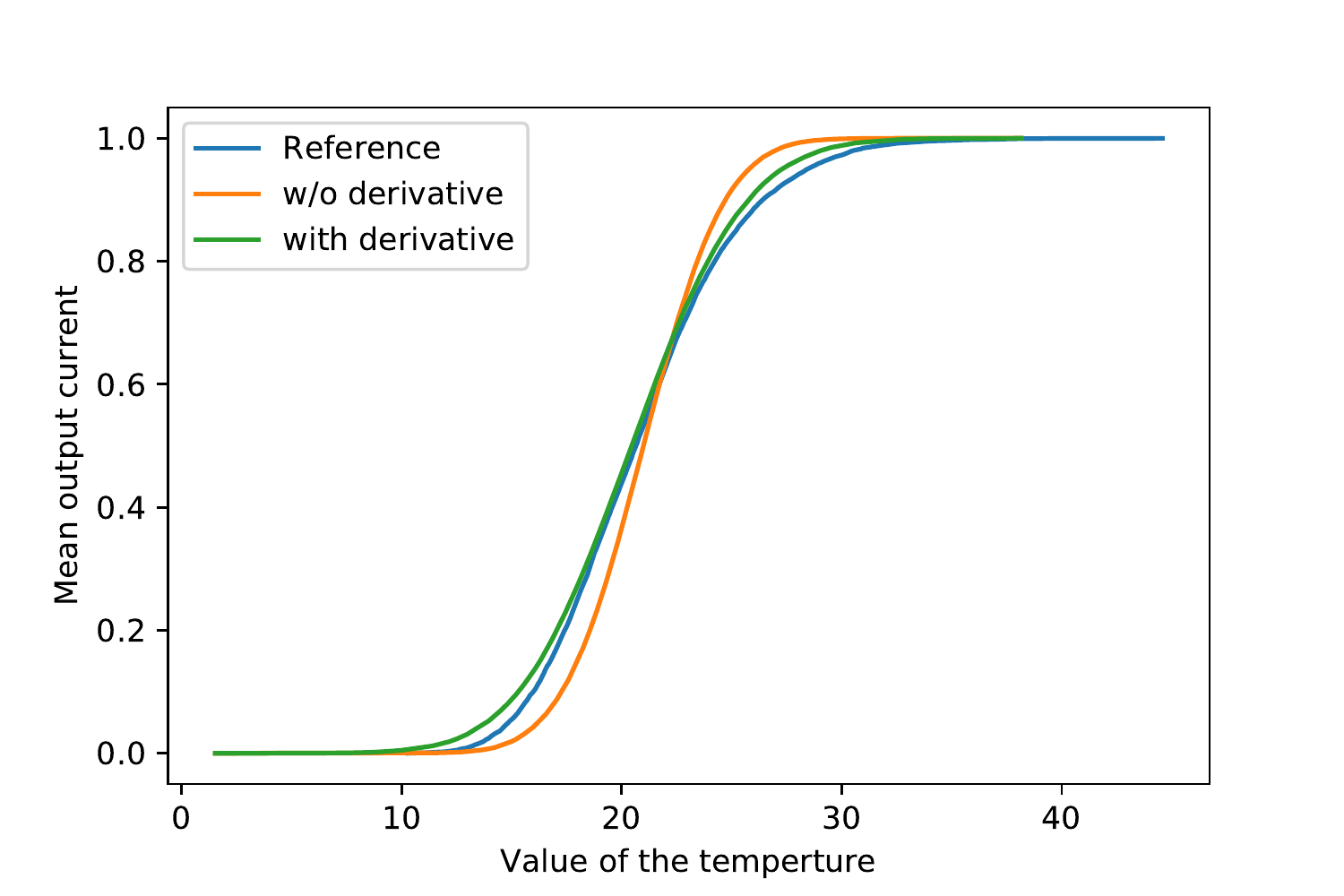}%
\caption{Cumulative distribution function of total temperature in the selected area.}%
\label{fig:RF_res}%
\end{figure}

The mesh used in the experiment is shown in Fig~\ref{fig:RF_mesh}.
We use Monte-Carlo simulation with $N_{MC}=10^4$ samples as a reference.
The results of numerical experiments are shown in Fig~\ref{fig:RF_res}.

The sizes of the system matrices~\eqref{system:mat} in cases of not using derivatives
and with derivatives are $20\times176$ and $81\times176\cdot41$, respectively.

We take $p=0.5$, and the total number of the basis polynomials is~$81$ when $q=2$.
In the non-derivative case, we just truncate polynomials with high powers.

\section{Possible difficulties of the derivative based approach}
\subsection{Singularity of system matrix}
In some cases, the system matrix is not of a full column rank.

Let us consider a small matrix generated by two-dimensional tensor product of  polynomials of  cumulative degree no more than~$2$ taken at two points.
It has the form
\begin{equation*}
    A=
    \left(
    \begin{array}{cccccc}
        1 & x_1 & {y_1} & {x_1} {y_1} & x_1^2 & y_1^2 \\
        0 & 1 & 0 & {y_1} & 2 {x_1} & 0 \\
        0 & 0 & 1 & {x_1} & 0 & 2 {y_1} \\
        1 & x_2 & {y_2} & {x_2} {y_2} & x_2^2 & y_2^2 \\
        0 & 1 & 0 & {y_2} & 2 {x_2} & 0 \\
        0 & 0 & 1 & {x_2} & 0 & 2 {y_2} \\
    \end{array}
    \right).
\end{equation*}
The rank of the matrix~$A$ is $5$ as the rows~$a_i$ for $i=1,\,\dots,6$ of the matrix are linearly dependent
\begin{equation*}
    - 2 a_1 + (x_1-x_2)a_2 + (y_1-y_2)a_3 +
    2 a_4 + (x_1-x_2)a_5 + (y_1-y_2)a_6 = 0.
\end{equation*}

A more complex example is the matrix of size~$10\times12$
generated by polynomials of cumulative degree no more than~$2$
taken at three points
\begin{equation*}
A=    \left(
\begin{array}{cccccccccc}
 1 & x_1 & y_1 & z_1 & x_1 y_1 & x_1 z_1 &
   y_1 z_1 & x_1^2 & y_1^2 & z_1^2 \\
 0 & 1 & 0 & 0 & y_1 & z_1 & 0 & 2 x_1 & 0 & 0 \\
 0 & 0 & 1 & 0 & x_1 & 0 & z_1 & 0 & 2 y_1 & 0 \\
 0 & 0 & 0 & 1 & 0 & x_1 & y_1 & 0 & 0 & 2 z_1 \\
 1 & x_2 & y_2 & z_2 & x_2 y_2 & x_2 z_2 &
   y_2 z_2 & x_2^2 & y_2^2 & z_2^2 \\
 0 & 1 & 0 & 0 & y_2 & z_2 & 0 & 2 x_2 & 0 & 0 \\
 0 & 0 & 1 & 0 & x_2 & 0 & z_2 & 0 & 2 y_2 & 0 \\
 0 & 0 & 0 & 1 & 0 & x_2 & y_2 & 0 & 0 & 2 z_2 \\
 1 & x_3 & y_3 & z_3 & x_3 y_3 & x_3 z_3 &
   y_3 z_3 & x_3^2 & y_3^2 & z_3^2 \\
 0 & 1 & 0 & 0 & y_3 & z_3 & 0 & 2 x_3 & 0 & 0 \\
 0 & 0 & 1 & 0 & x_3 & 0 & z_3 & 0 & 2 y_3 & 0 \\
 0 & 0 & 0 & 1 & 0 & x_3 & y_3 & 0 & 0 & 2 z_3 \\
\end{array}
\right)
\end{equation*}
The rank of this matrix~$A$ is~$9$.

This effect is also manifested in the matrices based on the polynomials of the maximal degree not more than some fixed value.
Let~$A$ be the $15\times15$ matrix of the polynomials with the maximal power no more than~$4$ taken at~$5$ points.
The rank of this matrix is~$14$.

Thus, it limits from below the number of points we take in the derivative-based approach.
At the moment we do not have explicit values of the rank
and plan to study it in the future.
\subsection{When matrix without derivatives works better}
Let us consider a toy example of 1D function approximation at two points.
Let~$f(x)=\sin(x)$, evaluation points are $x_0=-\pi$ and $x_1=\pi$.
We use a linear function $\widehat f(x)=kx+b$ to approximate~$f$.
The standard approach leads to the following linear system
\begin{equation*}
   \begin{pmatrix}
       1&-\pi\\
       1& \pi
   \end{pmatrix}
   \begin{pmatrix}
       b\\
       k
   \end{pmatrix}
=
   \begin{pmatrix}
       0\\
       0
   \end{pmatrix}
\end{equation*}
with solution~$b=0$, $k=0$.
The norms of the error of this approximation are
\begin{equation*}
    \norm|f-\widehat f|_2=\int_0^{2\pi}\sin^2 x\,dx\approx1.77,
\quad
    \norm|f-\widehat f|_\infty=1.
\end{equation*}
The linear system of the LSM in derivative approach is the following
\begin{equation*}
   \begin{pmatrix}
       1&-\pi\\
       0&1\\
       1& \pi\\
       0&1
   \end{pmatrix}
   \begin{pmatrix}
       b\\
       k
   \end{pmatrix}
=
   \begin{pmatrix}
       0\\
       -1\\
       0\\
       -1
   \end{pmatrix}
\end{equation*}
with the solution
\begin{equation*}
    b=0,
    \quad
    k=-\frac1{1+\pi^2}.
\end{equation*}
In this approach the norms of the errors are
\begin{equation*}
    \norm|f-\widehat f|_2\approx2.11,\\
\quad
    \norm|f-\widehat f|_\infty\approx1.15.
\end{equation*}
One can see that both norms are worse in the derivative approach.

However, we can approximate our function with more than two basis polynomial in derivatives approach using the same number of points.
Let $\widehat f(x) = a_0 + a_1 x + a_2 x^2 + a_3 x^3$.
Then, after solving the linear system, we get
\begin{equation*}
    a_0 = 0,\quad a_1 = \frac12,\quad a_2 = 0,\quad a_3 = -\frac1{2\pi^2},
\end{equation*}
and the corresponding errors are
\begin{equation*}
    \norm|f-\widehat f|_2\approx0.71,\\
\quad
    \norm|f-\widehat f|_\infty\approx0.43.
\end{equation*}

So, the advantages of the derivative approach are manifested in the case when we take an approximation polynomial of a higher degree than in the case without derivatives.
Moreover, in the derivative approach we can generate a system matrix with a large ratio of the number of rows to the number of columns to obtain greater stability.

\section{Related works}
Our method
is applicable
only in the cases where the function can be approximated by a linear combination of fixed basis functions.
The problem of function approximation using
a preselected set of basis functions
is the old and classical problem.
For many classes of problems,
including
the approximation of
smooth multivariate functions~\cite{T82,Z86},
solutions of parametric PDEs~\cite{CCS14},
and uncertainty quantification~\cite{BS11,CDS11},
approximation bounds are known.
The interpolation and approximation of functions in this setting using adaptive function sampling have been studied in~\cite{BH15},
including bounds on the Lebesgue constant for different basis sets in the multidimensional cases~\cite{GT14,CP11}.

There are several approaches to improve the accuracy of the approximation and to
obtain
theoretical estimates on it.
In the paper~\cite{D82} an approximation by smooth spline is considered in the multidimensional case.
An adaptive algorithm is proposed which allows selecting domains of approximation.
Errors estimations are given for this adaptive scheme.

The paper~\cite{CDS11} considers the regularity of the polynomial approximation to the solution of parametric PDEs.

An approximation of the solutions of the parametric PDEs is also considered in the paper~\cite{CD15}.
The solutions as functions of the parameters are holomorphic and highly anisotropic.
The convergence rate of the approximation is established using these properties of the solution,
which allows to avoid the curse of dimensionality.

A family of functions which are used for the approximation may not form an orthonormal basis.
In connection with this, a problem of finding appropriate basis occur.
In the paper~\cite{BCDD08} the approximation of functions in Hilbert space is considered.
A greedy algorithm that allows selecting of the basis functions is proposed.
Applications to the regression problem are considered.
For a class of greedy algorithms convergence results are proven.

The classical paper~\cite{S63} considers the interpolation of certain classes of multivariate functions.
An efficient numerical algorithm is developed based on this study.

One of the ways to increase the accuracy of the approximation and reduce the number of points of function evaluation
is finding appropriate sample points.
The paper~\cite{NJ2014} presents an adaptive sparse grid stochastic collocation approach.
Sparse grids interpolants are used to approximate solution as a function of input parameters.
Univariate Leja sequence is suggested to use for a high-dimensional interpolation.

The classical problem of reconstruction of an unknown function
by sampling in random points, and
with the use of the least square method is considered in the paper~\cite{CDL2013}.
A criterion on the number of approximating basis functions is presented.
This criterion allows to ensure the stability of the LSM,
and that the accuracy of the approximation is comparable to the best approximation.
This criterion can be applied for different approximation schemes including trigonometric and algebraic polynomials approximation.

The paper~\cite{BPS17} presents an algorithm for choosing points for the function interpolation.
This method is applied in the adaptive fashion to different functions
which are high-dimensional and
have high evaluation cost.

The paper~\cite{GI17} investigates a weighted method of least squares in the application to the function approximation problem
by truncating the polynomial chaos expansion.
The weights are calculated on the basis of the selected points where the function is evaluated.
These values correspond to the quadratures of integrals arising in the coefficient estimation.
Such an approach leads to the stability and high accuracy of the approximation.

The integration problem of a multivariate function is solved in the paper~\cite{SW98}
by using quasi-Monte Carlo (QMC) method.
The authors identify classes of functions for which the effect of dimension is negligible.
It is proved that the error in the worst case does not depend on the number of dimensions under certain conditions.

In the paper~\cite{KS05} several methods of avoiding the curse of dimensionality are considered in an application to the problem of integration over the $d$-dimensional cube.

The problems of multivariate function approximation and integration are considered in the paper~\cite{HW00}.
It is proved that the studied integration schemes are strongly tractable for a certain set of weighted function spaces.

In the survey~\cite{BG04} the theory and applications of sparse grids are presented.
The authors focus on the solution of partial differential equations.
Several estimations on the number of degrees of freedom are considered.
An extension to non-smooth solutions is made by adaptive update methods.

In our research, we rely on the fast methods for derivative computation
which have attracted tremendous attention in the recent years in the field of machine learning.
Many modern packages like TensorFlow~\cite{TF} and Pytorch~\cite{PyT} implement such tools.

\section{Conclusions and future work}
In the present paper advantages of multivariate function interpolation with the use of the information about its derivatives are demonstrated.
Several numerical experiments are conducted
including ODE and PDE equations with coefficients depending on multidimensional parameters.

As a subject of future work,
we consider the task of determining the optimal relation between the number of sample points
and the number of the basis functions
to increase stability and accuracy of the proposed scheme.
The problem of finding optimal points set of function evaluation
is in our plans as well.
Additionally,
we want to investigate convergence and stability bounds when the number of points increases.
Also, it is crucial to determine theoretical estimates for the Lebesgue constants for the derivative approach.

\bibliographystyle{siamplain}
\bibliography{refs}

\hyphenation{Post-Script Sprin-ger}
\begin{thebibliography}{10}

\bibitem{AHKAB}
{\em Ahkab python package}.
\newblock \url{https://ahkab.github.io/ahkab/}.
\newblock [online] [Accessed: 2017-03-20].

\bibitem{AG}
{\em Autograd python library}.
\newblock \url{https://github.com/HIPS/autograd/}.
\newblock [online] [Accessed: 2017-11-01].

\bibitem{DAdj}
{\em Dolfin-adjoin project}.
\newblock \url{http://www.dolfin-adjoint.org/en/latest/}.
\newblock [online] [Accessed: 2017-11-01].

\bibitem{PyT}
{\em Pytorch python framework}.
\newblock \url{https://www.tensorflow.org/}.
\newblock [online] [Accessed: 2017-11-01].

\bibitem{BM}
{\em ''rect maxvol`` algorithm}.
\newblock \url{https://bitbucket.org/muxas/rect_maxvol}.
\newblock [online] [Accessed: 2017-11-01].

\bibitem{TF}
{\em Tensorflow library}.
\newblock \url{https://www.tensorflow.org/}.
\newblock [online] [Accessed: 2017-11-01].

\bibitem{BCDD08}
{\sc A.~R. Barron, A.~Cohen, W.~Dahmen, and R.~A. DeVore}, {\em Approximation
  and learning by greedy algorithms}, Ann. Stat., 36 (2008), pp.~64--94.

\bibitem{BS83}
{\sc W.~Baur and V.~Strassen}, {\em The complexity of partial derivatives},
  Theor. Comput. Sci., 22 (1983), pp.~317--330.

\bibitem{BS11}
{\sc G.~Blatman and B.~Sudret}, {\em Adaptive sparse polynomial chaos expansion
  based on least angle regression}, Journal of Computational Physics, 230
  (2011), pp.~2345--2367.

\bibitem{BG04}
{\sc H.-J. Bungartz and G.~Michael}, {\em Sparse grids}, Acta Numer., 13
  (2004), pp.~1--123.

\bibitem{BPS17}
{\sc E.~Burnaev, I.~Panin, and B.~Sudret}, {\em Efficient design of experiments
  for sensitivity analysis based on polynomial chaos expansions}, Ann. Math.
  Artif. Intel., 81 (2017), pp.~187--207.

\bibitem{CP11}
{\sc J.-P. Calvi and V.~M. Phung}, {\em On the lebesgue constant of {Leja}
  sequences for the unit disk and its applications to multivariate
  interpolation}, J. Approx Theory, 163 (2011), pp.~608--622.

\bibitem{CCS14}
{\sc A.~Chkifa, A.~Cohen, and C.~Schwab}, {\em High-dimensional adaptive sparse
  polynomial interpolation and applications to parametric {PDEs}}, Found.
  Comput. Math., 14 (2014), pp.~601--633.

\bibitem{CDL2013}
{\sc A.~Cohen, M.~A. Davenport, and D.~Leviatan}, {\em On the stability and
  accuracy of least squares approximations}, Found. Comput. Math., 13 (2013),
  pp.~819--834.

\bibitem{CD15}
{\sc A.~Cohen and R.~DeVore}, {\em Approximation of high-dimensional parametric
  {PDEs}}, Acta Numer., 24 (2015), pp.~1--159.

\bibitem{CDS11}
{\sc A.~Cohen, R.~DeVore, and C.~Schwab}, {\em Analytic regularity and
  polynomial approximation of parametric and stochastic elliptic {PDE’s}},
  Anal. Appl., 9 (2011), pp.~11--47.

\bibitem{D82}
{\sc W.~A. Dahmen}, {\em Adaptive approximation by multivariate smooth
  splines}, J. Approx Theory, 36 (1982), pp.~119--140.

\bibitem{DTU16}
{\sc D.~Dung, N.~V. Temlyakov, and T.~Ullrich}, {\em Hyperbolic cross
  approximation},  (2016), \url{https://arxiv.org/abs/1601.03978}.

\bibitem{GI17}
{\sc S.~Ghili and G.~Iaccarino}, {\em Least squares approximation of polynomial
  chaos expansions with optimized grid points}, SIAM J. Sci. Comput., 39
  (2017), pp.~A1991--A2019.

\bibitem{G92}
{\sc A.~Griewank}, {\em Achieving logarithmic growth of temporal and spatial
  complexity in reverse automatic differentiation}, Optimization Methods and
  Software, 1 (1992), pp.~35--54.

\bibitem{GT14}
{\sc M.~Gunzburger and A.~L. Teckentrup}, {\em Optimal point sets for total
  degree polynomial interpolation in moderate dimensions},  (2014),
  \url{https://arxiv.org/abs/1407.3291}.

\bibitem{HW00}
{\sc F.~J. Hickernell and H.~Wozniakowski}, {\em Integration and approximation
  in arbitrary dimensions}, Adv. Comput. Math., 12 (2000), pp.~25--58.

\bibitem{KS05}
{\sc F.~Y. Kuo and I.~Sloan}, {\em Lifting the curse of dimensionality}, Not.
  Am. Math. Soc., 52 (2005), pp.~1320--1328.

\bibitem{EOLE}
{\sc C.-C. Li and A.~D. Kiureghian}, {\em Optimal discretization of random
  fields}, J. Eng. Mech., 119 (1993), pp.~1136--1154.

\bibitem{LHS}
{\sc M.~D. Mckay, R.~Beckman, and W.~Conover}, {\em A comparison of three
  methods for selecting vales of input variables in the analysis of output from
  a computer code}, Technometrics, 21 (1979), pp.~239--245.

\bibitem{MO16}
{\sc A.~Mikhalev and I.~V. Oseledets}, {\em Rectangular maximum-volume
  submatrices and their applications}, Linear Algebra Appl., 538 (2018),
  pp.~187--211.

\bibitem{NJ2014}
{\sc A.~Narayan and J.~D. Jakeman}, {\em Adaptive {Leja} sparse grid
  constructions for stochastic collocation and high-dimensional approximation},
  SIAM J. Sci. Comput., 36 (2014), pp.~A2952--A2983.

\bibitem{SW98}
{\sc I.~H. Sloan and H.~Wozniakowski}, {\em When are quasi-monte carlo
  algorithms efficient for high dimensional integrals?}, J. Complexity, 14
  (1998), pp.~1--3.

\bibitem{S63}
{\sc S.~A. Smolyak}, {\em Quadrature and interpolation formulas for tensor
  products of certain classes of functions}, Dokl. Akad. Nauk SSSR, 148 (1963),
  pp.~1042--1045.

\bibitem{T82}
{\sc V.~N. Temljakov}, {\em Approximation of periodic functions of several
  variables with bounded mixed difference}, Math. USSR Sb., 41 (1982),
  pp.~53--66.

\bibitem{BH15}
{\sc M.~Van~Barel and M.~Humet}, {\em Good point sets and corresponding weights
  for bivariate discrete least squares approximation}, Dolomites Research Notes
  on Approximation, 8 (2015), pp.~37--50.

\bibitem{Z13}
{\sc Z.~Zhang, T.~A. El-Moselhy, I.~M. Elfadel, and L.~Daniel}, {\em Stochastic
  testing method for transistor-level uncertainty quantification based on
  generalized polynomial chaos}, IEEE T. Comput. Aid. D., 32 (2013),
  pp.~1533--1545.

\bibitem{Z14}
{\sc Z.~Zhang, X.~Yang, I.~V. Oseledets, G.~E. Karniadakis, and L.~Daniel},
  {\em Enabling high-dimensional hierarchical uncertainty quantification by
  {ANOVA} and tensor-train decomposition}, IEEE T. Comput. Aid. D., 34 (2014),
  pp.~63--76.

\bibitem{Z86}
{\sc D.~Zung}, {\em Approximation of classes of smooth functions of several
  variables}, Journal of Soviet Mathematics, 35 (1986), pp.~2859--2875.

\end{thebibliography}

\end{document}